\documentclass[a4paper,11pt]{amsart}
\usepackage{amsmath,amsthm,amssymb}
\theoremstyle{plain}
\newtheorem{theorem}{Theorem}
\newcommand{\bem}{\begin{bmatrix}}
\newcommand{\eem}{\end{bmatrix}}
\DeclareMathOperator{\tr}{tr}
\DeclareMathOperator{\Id}{Id}
\DeclareMathOperator{\Rea}{Re}
\DeclareMathOperator{\Ima}{Im}
\newcommand{\marginextend}[1]{ \addtolength{\oddsidemargin}{-#1}  \addtolength{\evensidemargin}{-#1}
  \addtolength{\textwidth}{#1}\addtolength{\textwidth}{#1}}
\newcommand{\updownextend}[1]{ \addtolength{\topmargin}{-#1}  \addtolength{\textheight}{#1}
\addtolength{\textheight}{#1}}
\marginextend{1cm}
\updownextend{0cm}
\allowdisplaybreaks[4]
\title{A short proof of the elliptical range theorem}
\author{Gyula Lakos}
\address{Department of Geometry, Institute of Mathematics, E\"otv\"os University, P\'azm\'any P\'eter s.~1/C,  Budapest, H--1117, Hungary}
\email{lakos@cs.elte.hu}
\begin{document}
\begin{abstract}
A short proof of the elliptical range theorem concerning the numerical range of $2\times2$ complex matrices is given.
\end{abstract}
\maketitle
The basic properties of the numerical range were uncovered by Toeplitz \cite{T} and Hausdorff \cite{H}.
A key feature is the elliptical range theorem of $2\times2$ complex matrices due to Toeplitz.
We present a proof in the spirit of Toeplitz and Davis \cite{D}, but more explicitly,
including the lengths of the axes (cf. Uhlig \cite{U}),
in presentation more like as of Li \cite{L}.
\begin{theorem} If $A$ is a $2\times2$ complex matrix, then its numerical range
\[\mathrm W(A)=\{\langle A\mathbf x,\mathbf x\rangle\,:\,\mathbf x\in\mathbb C^2,\, |\mathbf x|=1 \}\]
is a possibly degenerate elliptical disk on the complex plane, with the eigenvalues of $A$ as the foci,
and with $s^{\pm}(A)=\frac12\sqrt{ \tr\,\left(A-\frac{\tr A}2\Id\right)^*\left(A-\frac{\tr A}2\Id\right)
\pm \left|\tr\,\left(A-\frac{\tr A}2\Id\right)^2\right|}$
as the major and minor semi-axes.
\begin{proof}
Applying the transform $A\mapsto \mathrm e^{i\theta}A+v\Id$, $\theta\in\mathbb R$, $v\in\mathbb C$,
it transforms the range and the eigenvalues accordingly, while $s^\pm(A)$ are left invariant.
Also, conjugation by a unitary matrix leaves all these data invariant.
By this we can assume that
$A=\bem c&2b\\0&-c\eem$ such that $b,c\in[0,+\infty)$.
Then, for $\mathbf x=\bem z_1\\z_2\eem$, $|z_1|^2+|z_2|^2=1$, it yields
\begin{align}
\label{eq}
\bem \Rea \langle A\mathbf x,\mathbf x\rangle
\\
 \Ima \langle A\mathbf x,\mathbf x\rangle
\\
0
\eem&=
\bem b&&-c\\&b&\\&&0 \eem
\bem2\Rea(\bar z_1z_2) \\ 2\Ima(\bar z_1 z_2)\\ |z_2|^2-|z_1|^2\eem
\\
&=\underbrace{\bem \sqrt{b^2+c^2}&&\\&b&\\&&0\eem}_F
\underbrace{\bem \frac{b}{\sqrt{b^2+c^2}}&&\frac{-c}{\sqrt{b^2+c^2}}
\\&1&\\\frac{c}{\sqrt{b^2+c^2}}&&\frac{b}{\sqrt{b^2+c^2}}\eem}_R
\underbrace{\bem2\Rea(\bar z_1z_2) \\ 2\Ima(\bar z_1z_2)\\ |z_2|^2-|z_1|^2\eem}_S
\notag
\end{align}
(in the case of $b=c=0$, any orthogonal matrix can be chosen for $R$).

Taking all unit vectors $\mathbf x$, $S$ ranges over the unit sphere (the base of the Hopf fibration).
Applying $R$ leaves it invariant.
Applying $F$ independently dilates in the first and second coordinates, and totally contracts in the third one.
Thus \eqref{eq} (third coordinate omitted) ranges over possibly degenerate elliptical disk
of canonical position with major semi-axis $\sqrt{b^2+c^2}$ and minor semi-axis $b$.
Its foci are, then, $(\pm c,0)$.
Now, these data are according to the statement of the theorem.
\end{proof}
\end{theorem}

\end{document}